\documentclass[12pt]{amsart}
\usepackage{amsmath}
\usepackage{amsthm}
\usepackage{amscd}

\newcommand{\Z}{{\mathbb Z}}
\newcommand{\C}{{\mathbb C}}

\newcommand{\fg}{{\mathfrak g}}
\newcommand{\fh}{{\mathfrak  h}}

\newtheorem{theorem}{Theorem}
\newtheorem{lemma}{Lemma}
\newtheorem{corollary}{Corollary}
\newtheorem{proposition}[theorem]{Proposition}
\theoremstyle{definition}

\newtheorem{example}{Example}

\theoremstyle{remark}
\newtheorem{remark}{Remark}
\newtheorem*{acknowledgement}{Acknowledgement}
\begin{document}

\title[Tensor products of irreducible representations]{Unique
decomposition of tensor products\\
of irreducible representations of simple algebraic groups}

\author{C.~S.~Rajan}

\address{Tata Institute of Fundamental 
Research, Homi Bhabha Road, Bombay - 400 005, INDIA.}
\email{rajan@math.tifr.res.in}


\begin{abstract} 
We show that a tensor product of irreducible, finite dimensional
representations of a simple Lie algebra over a field of characteristic
zero, determines the individual constituents uniquely. This is
analogous to the uniqueness of prime factorisation of natural
numbers. 

\end{abstract}

\maketitle

\section{Introduction} 

\subsection{} 
Let $\fg$  be a simple Lie algebra over $\C$. The main aim of this
paper is to prove the following unique factorisation of tensor
products of irreducible,  finite dimensional representations of $\fg$: 

\begin{theorem}\label{main}
Let $\fg$  be a simple Lie algebra over $\C$. Let $V_1, \cdots, V_n$
and $W_1, \cdots, W_m$ be non-trivial, irreducible, finite dimensional
$\fg$-modules. Assume that there is an isomorphism of the tensor
products, 
\[  V_1\otimes \cdots \otimes V_n\simeq W_1\otimes \cdots \otimes
W_m,\]
as $\fg$-modules. Then $m=n$, and there is a permutation $\tau$ of the
set $\{1, \cdots, n\}$, such that 
\[ V_i\simeq W_{\tau(i)},\]
as $\fg$-modules. 
\end{theorem}
The particular case which motivated the above theorem is the following
corollary:
\begin{corollary}\label{adjoint}
Let $V, ~W$ be irreducible $\fg$-modules. Assume that 
\[ {\rm End}(V)\simeq {\rm End}(W),\]
as $\fg$-modules. Then $V$ is either isomorphic to $W$ or the dual
$\fg$-module $W^*$. 
\end{corollary}

When $\fg={\mathfrak{sl}}_2$, and the number of components is at most
two, the theorem follows by comparing the highest and lowest weights
that occur in the tensor product. However, this proof seems difficult
to generalize (see Subsection \ref{PRV}).  The first main step towards
a proof of the theorem, is to recast the hypothesis as an equality of
the corresponding products of characters of the individual
representations occurring in the tensor product. A pleasant,
arithmetical  proof for ${\mathfrak{sl}}_2$  (see Proposition
\ref{sl2}),  indicates that we are on a right route. The proof in the
general case depends on the fact that the Dynkin diagram of a simple
Lie algebra is connected, and   proceeds  by induction on the rank of
$\fg$, using the fact that any simple Lie algebra of rank $l$, has a
simple subalgebra of rank $l-1$.  We  analyze the restriction of the
Weyl character formula of $\fg$ to the centralizer  of the  simple
subalgebra, by expanding along the characters of the central
$\mathfrak{gl}_1$. We  compare the coefficients, which are characters
of the simple subalgebra, of the highest and the second highest
degrees occurring in the product. The highest degree term is again the
character corresponding to a tensor product of irreducible
representations. The second highest degee term is  a sum of the
products of irreducible characters. To understand this sum, we again
argue by induction using character expansions. However, instead of
leading to further complicated sums, the induction argument
stabilizes, and we can formulate and prove a linear independence
property of products of characters of a particular type. Combining the
information obtained from the highest and the second highest degree
terms occurring in the product, we obtain the theorem.

The outline of this paper is as follows: first we recall some 
preliminaries about representations and characters of semisimple Lie
algebras. We then give  the proof for ${\mathfrak{sl}}_2$, and also of an
auxiliary result which comes up in the proof by induction. Although
not needed for the proof in the general case, we present the proof for
$GL_n$, since the ideas involved in the proof seem a bit more
natural. Here the numerator of the Weyl-Schur character formula appears as a 
determinant, which can be looked upon as a polynomial function on the
diagonal torus. The inductive argument arises upon expanding  this
function in one of the variables, the coefficients of which are given
by the numerators occurring in the Weyl-Schur character formula for
appropriate representations of $GL_{n-1}$. We then set up the
formalism for general simple $\fg$, so that we can carry over  the
proof for $GL_n$ to  the general case.

\begin{acknowledgement} I am  indebted to Shrawan Kumar for
many useful discussions during the early part of this work. I also
thank  S. Ilangovan, R. Parthasarathy,  D. Prasad, M. S. Raghunathan,
S. Ramanan and C. S. Seshadri for useful discussions. The arithmetical
application to Asai representations was suggested by D. Ramakrishnan's
work, who had proved a similar result for the usual degree two Asai
representations, and I thank him for conveying to me his results.
\end{acknowledgement}

\section{Preliminaries} 
We fix the notation and recall some of the relevant aspects of the
representation and structure theory of semisimple Lie algebras. We
refer to \cite{H}, \cite{S} for further details. 
\begin{enumerate}
\item Let $\fg$ be a complex semisimple Lie algebra, 
$\fh$ a Cartan subalgebra of
$\fg$, and  $\Phi\subset \fh^*$ the roots of the pair $(\fg, \fh)$. 
\item Denote by $\Phi^+\subset \Phi$, the subset of positive roots with respect
to some ordering of the root system, and by  
$\Delta$ a base for $\Phi^+$. 
\item Let $\Phi^*\subset \fh, ~\Phi^{*+}, ~\Delta^*$ be respectively 
the set of coroots,
positive coroots and fundamental coroots. Given a root
$\alpha\in \Phi$, $\alpha^*$ will denote the corresponding
coroot. 
\item Denote by $<.,.>:\fh\times \fh^*\to \C$ the duality pairing. For any root
$\alpha$,  we have $<\alpha^*, \alpha>=2$, and the pairing takes
values in integers when the arguments consist of roots and co-roots.  
\item  Given a root $\alpha$, by the properties of the root system, 
there are  reflections $s_{\alpha}, ~s_{\alpha^*}$ of $\fh^*, ~\fh$
respectively, defined by 
\[s_{\alpha}(u)=u-<\alpha^*, u>\alpha\quad  \text{and}\quad 
s_{\alpha^*}(x)=x-<x,\alpha>\alpha^*,\]
where $x\in \fh$ and $u\in \fh^*$. We have $s_{\alpha}(\Phi)\subset
\Phi$ and $s_{\alpha^*}(\Phi^*)\subset
\Phi^*$. 
\item Let $W$ denote the Weyl group of the root system. 
   The Weyl group $W$ is generated by the
reflections $s_{\alpha}$ for $\alpha\in \Delta$, subject to the
relations (see \cite[Theorem 2.4.3]{C})
\begin{equation}\label{genrel}
s_{\alpha}^2=1\quad \text{and}\quad
s_{\alpha}s_{\beta}s_{\alpha}=s_{s_{\alpha}(\beta)}, \quad \forall
~\alpha, ~\beta \in \Phi.
\end{equation}
In particular $s_{\alpha}$ and $s_{\beta}$ commute if
$s_{\alpha}(\beta)=\beta$. 
 There is a natural isomorphism between the Weyl groups of
the root system and the dual root system, given by $\alpha\mapsto
\alpha^*$ and $s_{\alpha}=^t\!s_{\alpha^*}$ the transpose of
$s_{\alpha^*}$. We identify the two actions of the Weyl group. 

\item Denote by $P\subset \fh^*$ the lattice of integral weights, given by 
\[ P=\{\mu\in \fh^*\mid \mu(\alpha^*)\in \Z, \forall \alpha\in
\Phi^*\}.\]
Dually we have a definition of the lattice of integral co-weights
$P^*$. 

\item Let  $P_+$ be  
the set of dominant, integral weights with respect to the
chosen ordering,   defined by 
\[ P_+= \{\lambda \in P\mid \lambda (\alpha^*)\geq 0, \forall \alpha\in
\Phi^{*+}\}.\]
 The
irreducible $\fg$-modules are indexed by elements in $P_+$, given by
highest weight theory. To each dominant, integral weight $\lambda$, we
denote the corresponding irreducible $\fg$-module with highest weight
$\lambda$ by $V_{\lambda}$. 
Let $l=|\Delta|$ be the rank of $\fg$. Index the collection of
fundamental roots  by $\alpha_1, \cdots, \alpha_l$. 
Denote by $\omega_1,  \cdots,
\omega_l$ (resp. $\omega_1^*,  \cdots,
\omega_l^*$), the set of fundamental weights (resp. fundamental
co-weights) defined by 
\[ \omega_i(\alpha_j^*)=\delta_{ij}\quad  \text{and}
\quad \omega_i^*(\alpha_j)=\delta_{ij}, \quad 1\leq i, j\leq l.\]
The fundamental weights form a $\Z$-basis for $P$. 

\item Let $l(w)$ denote the length of an element in the Weyl group,
given by the least length of a word in the $s_{\alpha}, ~\alpha\in
\Delta$ defining $w$. Let $\epsilon(w)=(-1)^{l(w)}$ be the sign
character of $W$. 

\item {\em Weyl character formula.} 
All the representations considered will be finite dimensional. 
Let $V$ be a $\fg$-module. With
respect to the action of $\fh$, we have a decomposition, 
\[ V=\oplus_{\pi\in \fh^*} V^{\pi},\]
\[ \text{where} \quad V^{\pi}=\{v\in V\mid hv=\pi(h)v, ~h\in \fh\}.\]
The linear forms $\pi$ for which $V^{\pi}$ are non-zero belong to the
weight lattice $P$, and these are the weights of $V$. Let $\Z[P]$
denote the group algebra of $P$, with basis indexed by $e^{\pi}$ for
$\pi\in P$. The (formal) character $\chi_V\in \Z[P]$ of $V$ is defined by,
\[ \chi_{{V}}= \sum_{\pi\in P}m(\pi)e^{\pi},\]
where $m(\pi)=\text{dim}(V^{\pi})$ is the multiplicity of $\pi$. 
The character is a ring homomorphism from the Grothendieck ring $K[\fg]$
defined by the representations of $\fg$ to the group algebra
$\Z[P]$. In particular, 
\[ \chi_{{V}\otimes {V'}}=\chi_{{V}}\chi_{{V'}}.\]
 The irreducible $\fg$-modules are indexed by elements in $P_+$, given by
highest weight theory. To each dominant, integral weight $\lambda$, we
denote the corresponding irreducible $\fg$-module with highest weight
$\lambda$ by $V_{\lambda}$, and the corresponding character by
$\chi_{{\lambda}}$.  Let 
\[ \rho=\frac{1}{2}\sum_{\alpha\in \Phi^+}\alpha=\omega_1+\cdots
 +\omega_l.\]
Define the Weyl denominator $D$ as,
\[ D=\sum_{w\in W}\epsilon(w) e^{w\rho} \in \Z[P].\]
The Weyl character formula for $V_{\lambda}$ is given by, 
\[ \chi_{\lambda}=\frac{1}{D}\sum_{w\in W}\epsilon(w)
e^{w(\lambda+\rho)}.\]
Let $S_{\lambda}=\sum_{w\in W}\epsilon(w)
e^{w(\lambda+\rho)}$ denote the numerator occurring in the Weyl
character formula.  We have $D=S_0$. 
\end{enumerate} 

We now recast the Main theorem. From the theory of characters, 
the Main theorem is equivalent to the following theorem: 
\begin{theorem} Let $\fg$ be a simple Lie algebra over $\C$. Assume
that there are positive integers 
$n\geq m$, and non-zero dominant
weights $\lambda_1, \cdots, \lambda_n, ~\mu_1, \cdots, \mu_m$  in
$P_+$ satisfying,
\begin{equation} \label{chareq0}
 S_{\lambda_1}...S_{\lambda_n}=S_{\mu_1}..S_{\mu_m}(S_0)^{n-m}. 
\end{equation}
Then $m=n$, and there is a permutation $\tau$ of the
set $\{1, \cdots, n\}$, such that 
\[ \lambda_i= \mu_{\tau(i)}, \quad 1\leq i\leq n.\]
\end{theorem}
We adopt a slight change in the notation. Assume $n\geq m$. Then
(\ref{chareq0}) can be rewritten as, 
\begin{equation}\label{chareq}
 S_{\lambda_1}...S_{\lambda_n}=S_{\mu_1}..S_{\mu_n},
\end{equation}
 where $\mu_i=0$ for $m+1\leq i\leq n$.

\subsection{$\mathfrak{sl}_2$ and PRV-components.}\label{PRV}
 Let $\fg= \mathfrak{sl}_2$. Let $V_n$ denote the irreducible
representation of $\mathfrak{sl}_2$ of dimension $n+1$, isomorphic to
the symmetric $n^{th}$ power $S^n(V_1)$  of the standard
representation $V_1$.   Suppose we have an isomorphism of
$\mathfrak{sl}_2$-modules,
\[ V_{n_1}\otimes V_{n_2}\simeq V_{m_1}\otimes V_{m_2}.\]
For any pair of positive integers $l\geq k$, we have the decompostion,
\[ V_k\otimes V_l\simeq  V_{l+k} \oplus V_{l+k-2}\oplus 
\cdots \oplus  V_{l-k}. \]
It follows that $n_1+n_2=m_1+m_2$ by comparing the highest weights.
Assuming $n_1\geq n_2$ and $m_1\geq m_2$, we have on comparing the
lowest weights occurring in the tensor product, that
$n_1-n_2=m_1-m_2$. Hence the theorem follows in this special case.

It is immediate from the hypothesis of the theorem,  that we have an
equality of the sum of the highest weights corresponding to the
irreducible modules $V_1, \cdots, V_n$ and $W_1,\cdots, W_m$
respectively.  The above proof for $\mathfrak{sl}_2$ suggests the use
of PRV-components: if $V_\lambda$ and $V_{\mu}$ are highest weight
finite dimensional $\fg$-modules  with highest weights $\lambda$ and
$\mu$ respectively, and $w$ is an element of the Weyl group, then it
is known that there is a Weyl group translate
$\overline{\lambda+w\mu}$ of the weight $\lambda+w\mu$, which is
dominant and such that the corresponding highest weight module
$V_{\overline{\lambda+w\mu}}$ is a direct summand  in the tensor
product module $V_\lambda\otimes V_{\mu}$ (see \cite{SK1}).  These are
the generalized Parthasarathy-Ranga Rao-Varadarajan
(PRV)-components. The standard PRV-component is obtained by taking
$w=w_0$, the longest element in the Weyl group. But the above proof
for $\mathfrak{sl}_2$ does not generalize, as the following example
for the simple Lie algebra $\mathfrak{sp}_6$ shows that it is not
enough to consider just the standard PRV-component:
\begin{example}
\[\fg =\mathfrak{sp}_6, \quad \fh=\C<e_1, e_2, e_3>, \quad \Delta=\{e_1-e_2,
e_2-e_3, 2e_3\}, \quad w_0=-1.\]
 Consider the following  highest weights on $sp_6$:
\begin{align*}
\lambda_1 & =6e_1+ 4e_2+2e_3 & \lambda_2 & =4e_1+2e_2 \\
\mu_1 & =6e_1+2e_2+2e_3 &  \mu_2 & =4e_1+4e_2 
\end{align*}
Clearly $\lambda_1+\lambda_2=\mu_1+\mu_2$. 
Since the Weyl group  contains sign changes, we see that
there exists an element of the Weyl group such that
$\lambda_1-\lambda_2=w(\mu_1-\mu_2)$. 
\end{example}
Thus we are led  to consider
generalized PRV-components. 
The problem with this approach 
is that although the standard PRV-component can be characterised as
the component on which the Casimir acts with the smallest eigenvalue, 
 there is no abstract 
 characterisation of the generalized PRV-component inside the tensor
product.  It is not clear that a generalized  PRV-component of one
side of the  tensor product,
is also  a PRV-component for the other tensor product. 
 Although the PRV-components occur with
`high' multiplicity \cite{SK2}, 
(greater than or equal to the order of the double coset
$W_{\lambda}\backslash W/W_{\mu}$, where $W_{\lambda}$ and $W_{\mu}$
are the isotropy subgroups of $\lambda $ and $\mu$ respectively), 
the converse is not
true. Even for $\mathfrak{sl}_2$, it does not seem easy to extend the above
proof when the number of components involved is more than two.

\section{GL(2)} 
The aim of this and the following 
section is to prove the main theorem in the context of
$GL(r)$: 
\begin{theorem} \label{glr}
Let $G=GL(r)$. 
Suppose $V\simeq V_{\lambda_1}\otimes \cdots \otimes V_{\lambda_n}$
and $W\simeq V_{\mu_1}\otimes \cdots \otimes V_{\mu_m}$ are tensor
products of irreducible representations with non-zero  highest weights
$\lambda_1, \cdots, \mu_m$. Assume that $V\simeq W$ as
$G$-modules. Then $n=m$ and there is a permutation $\tau$ of
$\{1,\cdots. n\}$ such that for $1\leq i\leq n$, 
\[ V_{\lambda_i}=V_{\mu_{\tau(i)}}\otimes {\rm det}^{\alpha_i},\]
for some integers $\alpha_i$. 
\end{theorem}
Upto twisting by a power of the
determinant, we can assume that the highest weight representations
$V_{\lambda}$  of $GL(r)$ are parametrised by their `normalized'
highest weights,
\[ \lambda =(a_1, \cdots a_r), ~~a_1\geq a_2 \geq \cdots\geq  a_r=0,\]
and $a_i$ are nonnegative integers.   It is enough then to show
under the hypothesis of the theorem, that the normalized highest
weights coincide.  
Let $x=(x_1,\cdots, x_r)$ be a
multivariable. We have the symmetric functions (Schur functions),
defined as the quotient of two determinants, 
\[ \chi_{\lambda}
(x)=\frac{|x_j^{a_i+r-i}|}{|x_j^{r-i}|}=\frac{S_{\lambda}}{D}, \]
where $S_{\lambda}$ denotes the determinant appearing in the numerator
and $D$ the standard Vandermonde determinant appearing in the
denominator.   
It is known that on the set of regular diagonal matrices the Schur
function $\chi_{\lambda}$ is equal to the character of $V_{\lambda}$.
Since we have assumed $a_n=0$, we have that the polynomials
$S_{\lambda}$ and $x_1$ are coprime, for any highest weight $\lambda$.

Hence by character theory,  the hypothesis of the
theorem can be recast as 
\begin{equation} 
 S_{\lambda_1}...S_{\lambda_n}=S_{\mu_1}..S_{\mu_n},
\end{equation}
 and where $\mu_i=0$ for
$m+1\leq i\leq n$.  

Write for $1\leq i\leq n$, 
\[
\begin{split}
\lambda_i& =(a_{i1}, a_{i2}, \cdots, a_{i(r-1)}, 0),\\
\mu_i& =(b_{i1}, b_{i2}, \cdots, b_{i(r-1)}, 0).
\end{split}
\]

\subsection{GL(2)} 

We present now the proof of the theorem for $GL(2)$.

\begin{proposition}\label{sl2}
Theorem \ref{glr} is true for $GL(2)$. 
\end{proposition}
\begin{proof}
Specializing Equation \ref{chareq} to the case of $GL(2)$, we obtain
\[ (x_1^{a_1+1}-x_2^{a_1+1})\cdots (x_1^{a_n+1}-x_2^{a_n+1})
=(x_1^{b_1+1}-x_2^{b_1+1})\cdots
(x_1^{b_n+1}-x_2^{b_n+1}), \]
where for the sake of simplicity we drop one of the indices in the
weights. Specialising the equation to $x_2=1$, and letting $x=x_1$, we
obtain an equality of the product of polynomials,
\[ (x^{a_1+1}-1)(x^{a_2+1}-1)\cdots (x^{a_n+1}-1)
=(x^{b_1+1}-1)(x^{b_2+1}-1)\cdots
(x^{b_n+1}-1). \]
Assume that $a_1={\rm max}\{a_1,\cdots, a_n\}$ and $b_1={\rm
max}\{b_1,\cdots, b_n \}$. For any positive integer $m$, let
$\zeta_m$ denote a primitive $m^{th}$ root of unity. 
The left hand side polynomial has a zero at $x=\zeta_{a_1+1}$, and the
equality forces the right side polynomial to vanish at
$\zeta_{a_1+1}$. Hence we obtain that $a_1\leq b_1$, and by symmetry
$b_1\leq a_1$. Thus $a_1=b_1$ and 
$\chi_{\lambda_1}=\chi_{\mu_1}$. Cancelling the first factor from both
sides, we are left with an equality of a product of characters
involving less number of factors than the equation we started with,
and by induction we have proved the theorem for $GL(2)$. 
\end{proof}

\begin{remark} It would be interesting to know the arithmetical
properties of the varieties defined by the polynomials $S_{\lambda}$
for general semisimple Lie algebras $\fg$. It seems difficult to
generalize the above arithmetical proof to general simple Lie
algebras. The proof in the general case proceeds by induction on the
rank, finally reducing  to the case of $\mathfrak{sl}_2$. 
\end{remark}

\subsection{A linear independence result} 
We now prove an auxiliary result for $GL(2)$, which arises in
the  inductive proof of Theorem \ref{glr}. 
\begin{lemma}\label{linindgl2} Let $\lambda_1, \cdots, \lambda_n$ be a
set of normalized weights in $P_+$.  Let $c$ be a positive integer and
$\omega_1$ denote the fundamental weight. Then the set
\[\{S_{\lambda_1}\cdots  S_{\lambda_{i-1}}
S_{\lambda_i+c\omega_1} S_{\lambda_{i+1}} \cdots   S_{\lambda_n}\mid
1\leq i\leq n\},\]
is linearly independent. In particular, suppose that there are 
 subsets $I, ~J\subset \{1, \cdots, n\}$
satisfying the following:
\begin{equation}
\begin{split} \label{multisumgl2}
 \sum_{i\in I} &  S_{\lambda_1}\cdots S_{\lambda_{i-1}}
S_{\lambda_i+c\omega_1} S_{\lambda_{i+1}}\cdots  S_{\lambda_n}\\
& = \sum_{j\in J}S_{\lambda_1}\cdots  S_{\lambda_{j-1}}
S_{\lambda_j+c\omega_1} S_{\lambda_{j+1}} \cdots S_{\lambda_n}.
\end{split}
\end{equation}
Then there is a bijection $\theta: I\to J$, such that
$\lambda_i=\lambda_{\theta(i)}$. 
\end{lemma}
 An equivalent statement can be made in the
Grothendieck ring $K[\fg]$ or with characters in place of
$S_{\lambda}$. 
\begin{proof}
Suppose we have a relation 
\[
\sum_{1\leq i\leq n} z_i S_{\lambda_1}\cdots  S_{\lambda_{i-1}}
S_{\lambda_i+c\omega_1} S_{\lambda_{i+1}} \cdots   S_{\lambda_n}
 =0,
\]
for some collection of complex numbers $z_i$. For any index $i$, let 
\[ E(i)=\{j\mid \lambda_j=\lambda_i\}.\] 
To show the linear independence, we have to show that for any index
$i$, we have 
\[ \sum_{j\in E(i)}z_j=0.\]
Dividing by $\prod_{l=1}^n S_{\lambda_l}$ on both sides and equating,
we are left with the equation,
\[ \sum_{1\leq i\leq n}z_i\frac{S_{\lambda_i+c\omega_1}}{S_{\lambda_i}}=0.\]
Specialising $x_1=1$ and writing $t$ instead of $x_2$, we obtain
\[ \sum_{1\leq i\leq n}z_i\frac{1-t^{a_i+c}}{1-t^{a_i}}=0.\]
Expand now as a power series in $t$. Consider the collection  of
indices $i$, for which $a_i$ attains the minimum value, say
for $i=1$. Equating the coefficient of $t^{a_1}$, we see  that 
$\sum_{j\in E(1)}
z_j=0$. Hence these terms can be removed from the relation, and we can
proceed by induction to  complete the proof of the lemma.
\end{proof}

\begin{remark} In retrospect, both Proposition \ref{sl2} and Lemma
\ref{linindgl2}, can be proved by comparing the coefficient of the
second highest power of $x_1$ occurring on both sides of the equation
(\ref{chareq}), as in the proofs occurring in the next section. But we
have included the proofs here, since it lays emphasis on the 
arithmetical properties  of the varieties defined by these
characters. 
\end{remark}

\section{Tensor products of $GL(r)$-modules}
We now come  to the proof of Theorem \ref{glr} for arbitrary $r$. 
The proof will proceed  by induction on $r$ and the maximum number of
components $n$, and we assume that the
theorem is true for $GL(s)$ with  $s<r$, and for $GL(r)$ with number
of components less than $n$.  
Associated to the highest weight $\lambda=(a_1, a_2, \cdots,
a_{r-1}, 0)$  of a  $GL(r)$-irreducible
module, define 
\[
\begin{split}
 \lambda' & =(a_2, a_3, \cdots, a_{r-1}, 0),\\
 \lambda'' & =(a_1+1, a_3, \cdots, a_{r-1}, 0)
\end{split}
\]
We can rewrite 
\[\lambda''=\lambda'+c(\lambda)\omega_1,\]
where  $\omega_1=(1,0, \cdots, 0)$ is 
the highest weight of the standard representation of $GL(r-1)$, and 
\begin{equation}
c(\lambda)=1+(a_1-a_2).
\end{equation} 
Both $\lambda'$ and $\lambda''$ are the highest weights of 
some $GL(r-1)$ irreducible modules. (Note: $0''\neq 0$). 

Expanding $S_{\lambda}$ as a polynomial in $x_1$ we obtain,
\[ \begin{split}
S_{\lambda}(x_1,\cdots, x_r) & = (-1)^{r+1}x_1^{a_1+r-1}
S_{\lambda'}(x_2,\cdots, x_r)\\
& +(-1)^r
x_1^{a_2+r-2}S_{\lambda''}(x_2,\cdots, x_r)+ Q,
\end{split}
\]
where $Q$ is a polynomial whose $x_1$ degree is less than $a_2+r-2$. 
Substituting in (\ref{chareq}),  and equating the
top degree term, we have an equality
\[
 x_1^{\sum_j a_{j1} +n(r-1)} S_{\lambda_1'}\cdots S_{\lambda_n'}=
 x_1^{\sum_j b_{j1} +n(r-1)}S_{\mu_1'}\cdots S_{\mu_n'}.
\]
 Since $S_{\lambda}$ and $x_1$ are
coprime polynomials for normalized $\lambda$, we have in particular
\begin{equation}\label{firstres}
S_{\lambda_1'}\cdots S_{\lambda_n'}= S_{\mu_1'}\cdots S_{\mu_n'}.
\end{equation}
Hence by induction it follows that there exists a permutation $\tau'$
of the set $\{1, \cdots, n\}$ such that,   
\[ \lambda_i'=\mu_{\tau'(i)}' \quad 1\leq i\leq n.\]

\begin{remark} At this point, with a little bit of extra work, a proof
of the main theorem can be given when the number of components $n$ is
at most two. Substituting $x_i=t^{i-1}$ for $1\leq i\leq n$, the
determinants $S_{\lambda}$ can be evaluated as Vandermonde
determinants. Arguing as in the proof of the main theorem for $GL(2)$,
it can be seen that we have an equality,
\[ \{a_{11}, a_{21}\}=\{b_{11}, b_{21}\}.\]
If we look at the constant coefficients, we obtain an equality,
\[ \{a_{11}-a_{12}, a_{21}-a_{22}\}=\{b_{11}-b_{12}, b_{21}-b_{22}\}.\]
These two inequalities combine to prove the main theorem when the
number of components is at most two. 
\end{remark}

\begin{remark} The proof does not proceed by restricting to various
possible $SL(2)$ mapping to $GL(r)$, and using the theorem for
$SL(2)$. For instance, it is not even possible to distinguish between
a representation and it's dual by restricting to various possible
$SL(2)$'s mapping to $GL(r)$. Morever we have the following example:   
\begin{example} Consider the following triples of highest weights on $GL(3)$:
\begin{align*}
\lambda_1 & =(3, 1,0) & \lambda_2 & =(2,2,0) &  \lambda_3&
=(1,0,0)\\
\mu_1 & =(3, 2,0) &  \mu_2 & =(2,0,0) &  \mu_3&
=(1,1,0)
\end{align*}
The calculations in the foregoing
remark, give in particular that for any co-root $\alpha^*$, the
sets of integers
\[\{<\alpha^*,\lambda_i+\rho>\mid 1\leq i\leq 3\}=
\{<\alpha^*,\mu_i+\rho>\mid 1\leq i \leq 3\} \] are equal, and it is
not  possible to differentiate the corresponding sets of highest
weights.  A calculation with the characters of the tensor product,
indicates that the term corresponding to the second highest
coefficient of $x_1$ in the corresponding  product of the characters
are different. This observation  motivates the rest of the proof of
Theorem \ref{main}.
\end{example}
\end{remark}
  
We continue with the proof of the theorem. Let
\[
 c_1 ={\rm min}\{c(\lambda_i)\mid 1\leq i \leq n\} \quad
\text{and}\quad  c_2  ={\rm min}\{c(\mu_j)\mid 1\leq j \leq n\}.
 \]
Equating  the coefficient of the second highest power of $x_1$ in
Equation \ref{chareq}, we obtain that $c_1=c_2=c$.  Let  $I, J\subset
\{1,\cdots, n\}$ denote the sets where for $i\in I$ and  $j\in J$,
$c(\lambda_i)$ (resp. $c(\mu_{\tau'(j)})$) attains the minimum
value. We obtain on equating the coefficient of the second highest
power of $x_1$ in Equation \ref{chareq}:
\[
 \sum_{i\in I}  S_{\lambda_1'}\cdots  S_{\lambda_{i-1}'}
S_{\lambda_i'+c\omega_2'} S_{\lambda_{i+1}'} \cdots   S_{\lambda_n'} =
\sum_{j\in J} S_{\lambda_1'} \cdots  S_{\lambda_{j-1}'}
S_{\lambda_j'+c\omega_2'} S_{\lambda_{j+1}'} \cdots S_{\lambda_n'}.
\]
Suppose there are indices $i\in I$ and $j\in J$ such that
$\lambda_i'=\mu_j'$ and $a_{i1}-a_{i2}=b_{j1}-b_{j2}=c-1$. It follows
that $\lambda_i=\mu_j$. Cancelling these two factors from the
hypothesis of Theorem \ref{glr}, we can proceed by induction on the
number of components to prove Theorem \ref{glr}. Thus, the proof of
Theorem \ref{glr}  reduces  now to the proof of the  following key
auxiliary lemma
(applied to  $GL(r-1)$), generalizing Lemma \ref{linindgl2}.
\begin{lemma} \label{linindglr}
Let $\lambda_1, \cdots, \lambda_n$ be a
set of normalized weights in $P_+$.  Let $d$ be a positive integer and
$\omega_1$ denote the fundamental weight. For any index $i$, let 
$E(i)=\{j\mid \lambda_j=\lambda_i\}$. 
Suppose we have a relation 
\[
\sum_{1\leq i\leq n} z_i S_{\lambda_1}\cdots  S_{\lambda_{i-1}}
S_{\lambda_i+d\omega_1} S_{\lambda_{i+1}} \cdots   S_{\lambda_n}
 =0,
\]
for some collection of complex numbers $z_i$. 
Then  for any index
$i$, we have 
\[ \sum_{j\in E(i)}z_i=0.\]
\end{lemma}

\begin{proof} By  Lemma \ref{linindgl2}, the lemma has been proved for
$GL(2)$. We argue by induction on the number of components $n$  and on $r$.
 Let
\[c={\rm min}\{c(\lambda_i)\mid 1\leq i \leq n\}.  \]
Let $M$ be the subset of $\{1, \cdots, n\}$ consisting of those
indices $i$ such that $c=1+ a_{i1}-a_{i2}$.  We expand both sides as a
polynomial in $x_1$, and compare the coefficient of the second highest
power of $x_1$. Since $d>0$, the term $S_{\lambda_i+d\omega_1}$ in the
product $S_{\lambda_1}\cdots  S_{\lambda_i+d\omega_1}\cdots
S_{\lambda_n}$ does not contribute to the coefficient of the second
highest degree term in $x_1$. Further we observe that
$(\lambda_i+d\omega_1)'=\lambda_i'$. Hence the contribution to the
coefficient of the second highest degree term in $x_1$ of
$S_{\lambda_1}\cdots  S_{\lambda_i+d\omega_1}\cdots S_{\lambda_n}$  is
given by,
\[\sum_{k\in M\backslash M\cap \{i\}} S_{\lambda_1'}\cdots 
S_{\lambda_k'+cf_1}\cdots 
S_{\lambda_n'},\]
where $f_1=(1,0,\cdots, 0)$ is a vector in $\Z^{r-1}$. Hence we
obtain, 
\[ \sum_{i=1}^nz_i  \sum_{k\in M\backslash M\cap\{i\}} S_{\lambda_1'}\cdots 
S_{\lambda_k'+cf_1}\cdots S_{\lambda_n'}=0. \]
Fix an index $i_0\in M$, and count the number of times
$\lambda_{i_0}'+cf_1$ terms occur as a component in the product. For
this, any index $i\neq i_0$ will contribute. Hence the above sum can
be rewritten as, 
\begin{equation}\label{auxind}
(n-1) \sum_{i\in M}z_i  S_{\lambda_1'}\cdots 
S_{\lambda_i'+cf_1}\cdots S_{\lambda_n'}=0. 
\end{equation} 
Cancelling those polynomials $S_{\lambda_l}'$ for $l\not\in M$, we
have by induction for any index $i_0\in M$, 
\[ \sum_{j\in E'(i_0)}z_j=0,\]
where $E'(i_0)=\{j\in M\mid \lambda_j'=\lambda_{i_0}'\}$.  We have
assumed that $j\in M$, since  we have cancelled those extraneous terms
with indices not in $M$. But since $j$ and $i_0$ belong to $M$, we
conclude that $\lambda_j=\lambda_{i_0}$. Hence $E'(i_0)=E(i_0)$, and
we have proved the lemma.
\end{proof}

\begin{remark} The surprising  fact is that the induction step rather
  than giving expressions where the character sums are spiked at more
  than one index, actually yields back Equation \ref{auxind}, which is
  again of the same type as the hypothesis of the lemma. 

 It is not clear whether there is a more
general context in which the above result can be placed. For example,
fix a weight $\Lambda$. Consider the  collection of characters 
\[\{\prod_{i\in I}\chi_{\lambda_i}\mid I\subset P_+, ~\sum_{i\in
I}\lambda_i=\Lambda\}.\]
It is not true that this set of characters is linearly independent,
since for a fixed $\Lambda$ the cardinality of this set grows
exponentially (since it is given by the  partition function), whereas
the dimension of the space of homogeneous polynomials in two variables
of fixed degree depends
polynomially (in fact linearly) on the degree.  It would be of
interest to generalize and place  the lemma in a more proper  context. 
\end{remark}

\section{Proof of the Main theorem in the general case} 
We now revert to the notation of Section 2. Our aim is to set up the
correct formalism  in the general case, so that we can carry over the
inductive proof for $GL(n)$ given above. Let $\fg$ be a simple Lie
algebra of rank greater than one.  Choose a fundamental root
$\alpha_1\in \Delta$.  Let $\Delta'=\Delta\backslash \{\alpha_1\}$,
and let $\Phi'\subset \Phi$ be the subset of roots lying in the span
of the roots generated by $\Delta'$. Let
\[ \fh'= \sum_{\alpha\in \Delta'}\C\alpha^*.\]
 It is known that $\Delta'$ is a base for the semisimple Lie algebra
$\fg'$ defined by,
\[\fg':=\fh'\oplus \sum_{\alpha\in
\Phi'}\fg^{\alpha}, \] where $\fg^{\alpha}$ is the weight space of
$\alpha$ corresponding to the adjoint action of $\fg$. The Lie algebra
$\fg'$ is a semisimple Lie algebra of rank $l-1$, and the roots of
$(\fg', \fh')$ can be identified with $\Phi'$. A special case is when
$\alpha_1$ corresponds to a corner vertex in the Dynkin diagram of
$\fg$. In this case,  $\fg'$ will be a simple Lie algebra.

We now want to find a suitable $\mathfrak{gl}_1$ complement to $\fg'$
inside $\fg$.  This is given by the following lemma:
\begin{lemma} 
$\omega_1^*\not\in \fh'$.
\end{lemma}
\begin{proof} Suppose we can write,
\[\omega_1^*=\sum_{i=2}^la_i\alpha_i^*,\]
for some complex numbers $a_i$. By definition of $\omega_1^*$, we
obtain on pairing with $\alpha_j $ for $2\leq j\leq l$,  the following
system of $l-1$ linear equations, in the unknown $a_i$:
\[   \sum_{i=2}^l<\alpha_i^*, \alpha_j>a_i=0.\]
But the matrix $(<\alpha_i^*, \alpha_j>)_{2\leq i, ~j\leq l}$ is the
Cartan matrix of the semisimple Lie algebra $\fg'$, and hence is
non-singular. Thus $a_i=0$ for $i=2,\cdots, l$, and that's a
contradiction as $\omega_1^*$ is non-zero.
\end{proof}

Let $W'$ denote the Weyl group of $(\fg',\fh')$, and it can be
identified with the subgroup of $W$ generated by the fundamental
reflections $s_{\alpha}$ for $\alpha\in \Delta'$. For such $\alpha$,
we have,
\[<s_{\alpha}(\omega_1^*),\alpha_j>=<\omega_1^*,\alpha_j>-<\omega_1^*,
\alpha><\alpha^*, \alpha_j>=\delta_{1j}.\] Hence  $W'$ fixes
$\omega_1^*$. Conversely, it follows from the fact that $\omega_1^*$
is orthogonal to all the roots $\alpha_2, \cdots, \alpha_l$, that any
element of the Weyl group fixing $\omega_1^*$ lies in the subgroup of
$W$ generated by the simple reflections $s_{\alpha_i}, ~2\leq i\leq
l$, and hence lies in $W'$ \cite[Lemma 2.5.3]{C}.

Our next step is to study the restriction of the character
$\chi_{\lambda}$ to $\fg'\oplus \mathfrak{gl}_1\subset \fg$.  Let $P'$
denote the lattice of weights of $\fg'$. We consider $P'$  as a
subgroup of $P$, consisting of those weights which vanish when
evaluated on $\omega_1^*$.

Choose a natural number $m$ such that $\pi(\omega_1^*)\in
\frac{1}{m}\Z$ for all weights $\pi\in P$. Let $Z_1$ be the subgroup,
isomorphic to the integers, of linear forms on $\fh$, which are
trivial on $\fh'$ and evaluated on $\omega_1^*$ lie in
$\frac{1}{m}\Z$. We have,
\[ P\subset Z_1\oplus P',\]
and we decompose the character $\chi_{\lambda}$ with respect to this
direct sum decomposition. Given a weight $\pi\in P$, denote by
$\pi'\in P'$ it's restriction to $\fh'$. Let $l_1$ be the weight in
$P$, vanishing on $\fh'$ and taking the value $1$ on $\omega_1^*$.
We define 
\[ d_1(\pi)=\pi(\omega_1^*)\]
as the {\em degree} of $\pi$ along $l_1$.
Write any weight $\pi$ with respect to the above decomposition as,
\[\pi=d_1(\pi)l_1 +\pi'\quad \text{so that} \quad
e^{\pi}=e^{d_1(\pi)l_1}e^{\pi'}. \]
The numerator of the Weyl character formula decomposes as,
\begin{equation}\label{charexpgen}
\begin{split}
 S_{\lambda}& =\sum_{d\in \frac{1}{m}\Z}e^{dl_1}\left(\sum_{w\in
W_d}\epsilon(w)e^{w(\lambda+\rho)'}\right),\\ \text{where}\quad  W_d&
=\{w\in W\mid (w(\lambda+\rho))(\omega_1^*)=d\}.
\end{split}
\end{equation}
We refer to the inner sum as the coefficient of the degree $d$
component along  $l_1$, or as the coefficient of $e^{dl_1}$.
Given a dominant integral
weight $\lambda\in P_+$, define
\[
\begin{split}
a_1(\lambda) & =\text{max}\{w\lambda(\omega_1^*)\mid w\in W\}, \\
 a_2(\lambda)&
=\text{max}\{w\lambda(\omega_1^*)\mid w\in W~ \text{and} ~
w\lambda(\omega_1^*)\neq a_1(\lambda)\}.
\end{split}
\]
The formalism that we
require in order to carry over the proof for $GL(n)$ to the general
case, is given by the following lemma:
\begin{lemma} \label{formalism}
Let $\lambda$ be a regular weight in $P_+$.  
\begin{enumerate}
\item The largest value
$a_1(\lambda)$ of $(w\lambda)(\omega_1^*)$ for $w\in W$,  is attained
precisely for 
$w\in W'$. In particular, 
\[ a_1(\lambda)=\lambda(\omega_1^*).\] 
\item The second highest value of $a_2(\lambda)$ is attained
precisely for $w$ in $W's_{\alpha_1}$, and the value is given by
\[a_2(\lambda)=s_{\alpha_1}\lambda(\omega_1^*)=
a_1(\lambda)-\lambda(\alpha_1^*).\] 
\end{enumerate}
\end{lemma}
\begin{proof}
1) By \cite[Lemma 2.5.3]{C}, we have to show that if $w\lambda(\omega_1^*)$
attains the maximum value, then $w$ fixes $\omega_1^*$. Since
$\omega_1^*$ is a fundamental co-weight, we have 
\[\omega_1^*-w\omega_1^*=\sum_{i=1}^ln_i\alpha_i^*, \]
for some non-negative natural numbers $n_i$. Hence, 
\[
\lambda(\omega_1^*-w\omega_1^*)=\sum_{i=1}^ln_i\lambda(\alpha_i^*),\]
and the latter expression is strictly positive, if some $n_i>0$, since
$\lambda$ is regular and dominant. This proves the first part.

2)  We prove the second part by induction on the length $l(w)$ of $w$. The
statement is clear for the fundamental reflections, which are of
length one. Consider an element of the form $ws_{\beta}$ such that 
$ws_{\beta}\lambda(\omega_1^*)=a_2(\lambda)$, $\beta$
is a fundamental root and 
$l(ws_{\beta})=l(w)+1$. By \cite[Lemma 2.2.1 and Theorem 2.2.2]{C},  
it follows that, 
\[ w(\beta) \in \Phi^+.\]
\[
\begin{split} \text{We have} \quad \lambda-ws_{\beta}\lambda &= 
(\lambda-w\lambda)+w(\lambda-s_{\beta}\lambda)\\
& =(\lambda-w\lambda)+<\beta^*, \lambda>w\beta.
\end{split}
\] 
Since $w\beta\in \Phi^+$, and $<\beta^*, \lambda>$ is positive, it
follows that 
\[ ws_{\beta}\lambda(\omega_1^*)\leq
w\lambda(\omega_1^*),\]
with strict inequality if $<\omega_1^*, w\beta>$ is positive. Hence by
induction we can assume that $w$ either belongs to $W'$,  or is of the
form $w_0s_{\alpha_1}$,  with $w_0\in W'$. We have to only consider the
second possiblity. We obtain, 
\[ (\lambda-w\lambda)(\omega_1^*)=
(\lambda-s_{\alpha_1}\lambda)(\omega_1^*)=<\alpha_1^*, \lambda>~ >0.\]
Assuming the hypothesis of the Lemma for $ws_{\beta}$, it follows that
$ws_{\beta}(\omega_1^*)=w(\omega_1^*)$, and hence   we
obtain  that $s_{\alpha_1}(\beta)$ has no $\alpha_1$ component,
when we expand it as a linear combination of the fundamental
roots. But
\[ s_{\alpha_1}(\beta)=\beta-<\alpha_1^*, \beta>\alpha_1,\]
and it follows that $<\alpha_1^*, \beta>=0$. From the relations
defining the Weyl group, it follows that $s_{\alpha_1}$ and
$s_{\beta}$ commute, and hence the element $ws_{\beta}$ is of the form
$w_1s_{\alpha_1}$ for some element $w_1\in W'$, and that concludes the
proof of the lemma. 
\end{proof}

The restriction of the fundamental weights $\omega_2,
\cdots, \omega_l$ to $\fh'$ are the fundamental weights (with the
indexing set ranging from $2, \cdots, l$, instead of $1, \cdots,
l-1$) of $\fg'$. In particular, the restriction $\rho'$ of $\rho$ is
the sum of the fundamental weights of $\fg'$. 
For $\lambda\in P_+$, define $\lambda''\in P_+'$ by,
\begin{equation}\label{lambda2or}
 \lambda''=(s_{\alpha_1}(\lambda+\rho))'-\rho'.
\end{equation}
We have the following corollary, giving the character expansion for
the first two terms along $l_1$:
\begin{corollary}\label{charexpdeg}
With notation as in 
the character expansion  given by Equation (\ref{charexpgen}), we have 
\begin{equation}
S_{\lambda} =e^{a_1(\lambda+\rho)l_1}S_{\lambda'}
- e^{a_2(\lambda+\rho)l_1}S_{\lambda''}+L(\lambda),
\end{equation}
where $L(\lambda)$ denotes the terms of degree along $l_1$ less than
the second highest degree. 
\end{corollary}
\begin{proof} The proof is immediate from Lemma \ref{formalism} and
equation \ref{charexpgen}, upon observing that $\lambda+\rho$ is
regular. The second term occurs with the opposite sign, since
$l(ws_{\alpha_1})=l(w)+1$, for $w\in W'$, and the length function of
$W$ restricts to the length function of $W'$, taken with respect to
$\Delta$ and $\Delta'$ respectively. 
\end{proof}

We write these down explicitly in terms of the fundamental weights. 
\[ \lambda=n_1(\lambda)\omega_1+\cdots+n_l(\lambda)\omega_l,\]
in terms of the fundamental weights, so that
$\lambda+\rho=(n_1(\lambda)+1)\omega_1+\cdots+(n_l(\lambda)+1)\omega_l$.
 We have, 
\begin{equation}
 (\lambda+\rho)' =(n_2(\lambda)+1)\omega_2'+\cdots+(n_l(\lambda)+1)\omega_l'.
\end{equation}
  
\[\text{Let}\quad  \fg'\simeq \oplus_{s\in S}\fg_s',\]
be the decomposition of $\fg'$ into simple Lie algebras. 
For each simple component $\fg'_s$ of $\fg'$,  let $\alpha_s$ be the
unique simple root connected to $\alpha_1$ in the Dynkin diagram of
$\fg$. Then 
\[- <\alpha_1^*, \alpha_s>=m_{1s},\]
is positive for each $s$.  This is possible, since we have
 assumed that $\fg$ is a
{\em simple} Lie algebra of larger rank! A calculation yields, 
\begin{equation}\label{lambda2}
s_{\alpha_1}(\lambda+\rho)'
=\sum_{s\in S}
(n_s(\lambda)+1+m_{1s}(n_1(\lambda)+1))\omega_s'+\sum_{t\in
\Delta\backslash S}n_t\omega_t'.
\end{equation}
For example, if $\alpha_1$ is a corner root, then $\fg'$ is
simple. Let $\alpha_2$ be the root adjacent to $\alpha_1$. In this
particular case, we have
\begin{equation}\label{lambda2corner}
 \lambda''=\lambda'+m_{12}(n_1(\lambda)+1)\omega_2'.
\end{equation}

We are now in a position to prove the main theorem, the proof of which
is along the same lines as the proof for $GL(n)$. We assume that we
have an equality as in equation (\ref{chareq}):
\begin{equation} 
 S_{\lambda_1}...S_{\lambda_n}=S_{\mu_1}..S_{\mu_n},
\end{equation} 
We now choose a corner root $\alpha_1$,  
and from  equation (\ref{charexpdeg}), we obtain,
\[\begin{split}
& \prod_{i=1}^n   (e^{a_1(\lambda_i+\rho)l_1}S_{\lambda_i'}-
e^{a_2(\lambda_i+\rho)l_1}S_{\lambda_i''}+L(\lambda_i))\\
& = \prod_{i=1}^n 
(e^{a_1(\mu_i+\rho)l_1}S_{\mu_i'}-e^{a_2(\mu_i+\rho)l_1}
S_{\mu_i''}+L(\mu_i)).
\end{split}
\]
 On taking
products and comparing the coefficients of the topmost degree, we
get,
\[  S_{\lambda_1'}...S_{\lambda_n'}=S_{\mu_1'}..S_{\mu_n'}.\]
By induction, we can thus assume that upto a permutation we have an
equality, 
\begin{equation} \label{inductiveeq}
(\lambda_1', \cdots, \lambda_n')=(\mu_1', \cdots, \mu_n').
\end{equation}
Now we compare the term contributing to the second highest degree in
the product. Let $I$ (resp. $J$) consist of those indices in the set
$\{1, \cdots, n\}$, for which $n_1(\lambda_i)=\lambda_i(\alpha_1^*)$ (resp.
$n_1(\mu_i)$) is minimum. By Part (2) of Lemma \ref{formalism} and the
character expansion as given by Corollary \ref{charexpdeg}, 
the second highest degree along $l_1$ 
in the product is of degree 
$n_1(\lambda+\rho)=n_1(\lambda)+1$ less than the total degree.  In
particular the minimum of $n_1(\lambda_i)$ and the minimum of
$n_1(\mu_j)$ coincide, as we vary over the indices. 
We have the following equality of the
second highest degree terms along $l_1$:
\[\sum_{i\in I}  S_{\lambda_1'}\cdots  S_{\lambda_{i-1}'}
S_{\lambda_i''} S_{\lambda_{i+1}'} \cdots   S_{\lambda_n'}
= \sum_{j\in J} S_{\mu_1'} \cdots  S_{\mu_{j-1}'}
S_{\mu_j''} S_{\mu_{j+1}'} \cdots S_{\mu_n'}.
\]
From equation (\ref{lambda2corner}), and equation (\ref{inductiveeq}), we can
recast this equality as,
\begin{equation}\label{recast}
\begin{split}
 & \sum_{i\in I} S_{\lambda_1'}\cdots S_{\lambda_{i-1}'}
S_{\lambda_i'+d\omega_2'}S_{\lambda_{i+1}'} \cdots   S_{\lambda_n'}\\
& = \sum_{j\in J} S_{\lambda_1'} \cdots S_{\lambda_{j-1}'}
S_{\lambda_j'+d\omega_2'}S_{\lambda_{j+1}'}\cdots S_{\lambda_n'}
\end{split}
\end{equation}
where $d=m_{12}(n_1(\lambda_i)+1)$ is a positive integer,  since
$\fg$ has been assumed to be simple.

Granting Lemma \ref{linindgeneral} given below,   the Main Theorem
follows, since we have indices $i_0, ~j_0$ such that,
\[ \lambda_{i_0}'=\mu_{j_0}' \quad \text{and}\quad
n_1(\lambda_{i_0})=n_1(\mu_{j_0}),\] where the indices $i_0, ~j_0$ are
such that the minimum of $n_1(\lambda_i)$ and $n_1(\mu_j)$ is
attained. Hence we have,
\[\lambda_{i_0}=\mu_{j_0}.\]
Cancelling these terms from equation (\ref{chareq}), we are left with
an equality where the number of components occurring in the tensor
product in the hypothesis of the Main theorem, is less than the one we
started with, and an induction on the number $n$ of components in the
tensor product proves the main theorem.

\begin{remark} When the number of components is at most two, then we
do not need Lemma \ref{linindgeneral}.  In  the above equality
(\ref{recast}), we can assume that $I\cap J$ is empty, and so can take
for example $I=\{1\}$ and $J=\{2\}$, to obtain,
\[
S_{\lambda_1'+d\omega_2'}S_{\lambda_2'}=S_{\lambda_1'}
S_{\lambda_2'+d\omega_2'}.
\]
By induction on the rank, assuming that the Main theorem is true with
number of components at most two, we obtain the Main theorem for all
$\fg$.
\end{remark}

To complete the proof of the theorem, we have to state the auxiliary
 linear independence property,   generalizing Lemma \ref{linindgl2}
 and Lemma \ref{linindglr}.

\begin{lemma} \label{linindgeneral} Let $\fg$ be a simple Lie algebra, and
let $\lambda_1, \cdots, \lambda_n$ be a set of dominant, integral
weights in $P_+$.   Let $d$ be a positive integer and $\omega_p$
denote a fundamental weight corresponding to the root $\alpha_p$.  For
any index $i$, let  $E(i)=\{j\mid \lambda_j=\lambda_i\}$.  Suppose we
have a relation
\[
\sum_{1\leq i\leq n} z_i S_{\lambda_1}\cdots  S_{\lambda_{i-1}}
S_{\lambda_i+d\omega_p} S_{\lambda_{i+1}} \cdots   S_{\lambda_n} =0,
\]
for some collection of complex numbers $z_i$.  Then  for any index
$i$, we have
\[ \sum_{j\in E(i)}z_i=0.\]
\end{lemma}
\begin{remark} Instead of $\omega_p$, we can spike up the equation
with any non-zero highest weight $\lambda$, but the proof is
essentially the same.
\end{remark}

The proof of this lemma will be by induction on the rank. For simple
Lie algebras not of type $D$ or $E$, and if $\omega_p$ is a
fundamental weight corresponding to a corner root in the Dynkin
diagram of $\fg$, the proof follows along the same lines as in the
proof of Lemma \ref{linindglr}, and that is sufficient to prove the
main theorem in these cases. For Lie algebras of type $D$ and $E$, the
proof becomes complicated, due to the fact that the root adjacent to
a corner root $\alpha_1$ in the Dynkin diagram of $\fg$, need not be a
corner root in the Dynkin diagram associated to $\Delta\backslash
\{\alpha_1\}$.  Before embarking on a proof of this lemma, we will
need a preliminary lemma. We say that a simple Lie algebra $\fg$ has
{\em Property LI}, if Lemma \ref{linindgeneral} holds for $\fg$.
\begin{lemma}\label{linindsemisimple}
Assume that Property LI holds for all simple Lie algebras of rank at
most $l$. Let $\oplus_{s\in S}\fg_s$ be a direct sum of simple Lie
algebras of $\fg_s$ of  rank at most $l$. For each $s\in S$, assume
that we are given dominant, integral weights $\lambda_{s1}, \cdots,
\lambda_{sn}$ of $\fg_s$, a positive integer $d_s$, and a fundamental
weight $\omega_s$  of $\fg_s$. Suppose that we have a relation,
\begin{equation}\label{semisimplemultisum}
\sum_{1\leq i\leq n} z_i S_{\Lambda_1}\cdots  S_{\Lambda_{i-1}}
S_{\hat{\Lambda}_i} S_{\Lambda_{i+1}} \cdots   S_{\Lambda_n} =0,
\end{equation}
for some collection of complex numbers $z_i$, where for $1\leq i\leq n$
\[ \begin{split}
S_{\Lambda_i} & =\prod_{s\in S}S_{\lambda_{si}}, \\ \text{and}\quad
S_{\hat{\Lambda}_i} & =\prod_{s\in S}S_{\lambda_{si}+d_s\omega_s}.
\end{split}
\]
Then  for any index $i$, we have
\[ \sum_{j\in E(i)}z_j=0,\]
where $E(i)=\{j\mid (\lambda_{sj})_{s\in S}= (\lambda_{si})_{s\in
S}\}$.
\end{lemma}
\begin{proof}
The proof proceeds by induction on the cardinality of $S$.  Consider
equation (\ref{semisimplemultisum}),  as a equation with respect to
one of the simple Lie algebras, say  $\fg_1$. The linear independence
property, reduces to the case when the number of simple Lie algebras
involved is one less, and we are through by induction.
\end{proof}

Now we get back to the proof of Lemma \ref{linindgeneral}.
\begin{proof}  By  Lemma \ref{linindgl2}, the lemma has been proved for
$GL(2)$. We argue by induction on the number of components $n$  and on
 the rank $l$ of $\fg$. We assume that the lemma has been proved for
 all simple Lie algebras of rank less than the rank of $\fg$. We use
 the character expansion given by Corollary \ref{charexpdeg}, where we
 denote by $l_p$ the linear functional corresponding to
 $\omega_p^*$. Let
\[ \fg'\simeq \oplus_{s\in S}\fg_s,\]
be the decomposition of $\fg'$ into simple Lie algebras (we are
removing the fundamental root $\alpha_p$ corresponding to the
fundamental weight $\omega_p$ from the Dynkin diagram). For each $s\in
S$, let $\alpha_s\in \Delta$ be the root adjacent to $\alpha_p$.

Let $M$ be the subset of $\{1, \cdots, n\}$ consisting of those
indices $i$ such that $n_p(\lambda_i)$ attains the  minimum.  For each
$s\in S$,    let
\[c_s={\rm min}\{m_{ps}(n_p(\lambda_i)+1)\mid i\in M\}.  \]
We expand both sides using Corollary (\ref{charexpdeg}),  and compare
the coefficient of the second highest degree along  $l_p$.  Since
$d>0$, the term $S_{\lambda_i+d\omega_p}$ in the product
$S_{\lambda_1}\cdots  S_{\lambda_i+d\omega_p}\cdots S_{\lambda_n}$
does not contribute to the coefficient of the second highest degree
term in $l_p$. Further we observe that
$(\lambda_i+d\omega_1)'=\lambda_i'$. Hence the contribution to the
coefficient of the second highest degree term in $l_p$ of
$S_{\lambda_1}\cdots  S_{\lambda_i+d\omega_p}\cdots S_{\lambda_n}$  is
given by,
\[\sum_{k\in M\backslash M\cap \{i\}} S_{\lambda_1'}\cdots 
S_{\lambda_{k-1}'}S_{\lambda_k''}S_{\lambda_{k+1}'}\cdots
S_{\lambda_n'}, \]
where $\lambda_k''$ is given by equation (\ref{lambda2}) as follows:
\[ \lambda_k''+\rho'=\sum_{s\in S}
(n_s(\lambda_k)+1+m_{ps}(n_p(\lambda_k)+1))\omega_s'+\sum_{t\in
\Delta\backslash S}n_t\omega_t'.
\]
Since $\fg$ is simple and $\lambda_k+\rho$ is regular, we notice that 
the term $m_{ps}(n_p(\lambda_k+\rho))$ is always positive.  On
rearranging the sum, we
obtain, 
\[ \sum_{i=1}^nz_i  \sum_{k\in M\backslash M\cap\{i\}} S_{\lambda_1'}\cdots 
S_{\lambda_k''}\cdots S_{\lambda_n'}=0. \]
Fix an index $i_0\in M$, and count the number of times a given 
index $\lambda_{i_0}''$ occurs as a component in the product. For
this, any index $i\neq i_0$ will contribute. Hence the above sum can
be rewritten as, 
\[(n-1) \sum_{i\in M}z_i  S_{\lambda_1'}\cdots 
S_{\lambda_i''}\cdots S_{\lambda_n'}=0. \] 
Cancel those polynomials $S_{\lambda_l}'$ for $l\not\in M$.  We
have by Lemma \ref{linindsemisimple}, for any index $i_0\in M$, 
\[ \sum_{j\in E'(i_0)}z_j=0,\]
where $E'(i_0)=\{j\in M\mid \lambda_j'=\lambda_{i_0}'\}$.
where we can assume that $j\in M$, since  we have cancelled those
extraneous terms with indices not in $M$. But since $j$ and $i_0$
belong to $M$ and $j$ is in $E(i_0)$,  
we conclude that $\lambda_j=\lambda_{i_0}$. Hence
$E'(i_0)=E(i_0)$, and we have proved the lemma.

\end{proof}

\begin{remark} The main theorem indicates the presence of an
`irreduciblity property' for the characters of irreducible
representations of simple algebraic groups. However the naive feeling
that the characters of irreducible representations are irreducible is
false. This can be seen easily for $\mathfrak{sl}_2$. For $GL(n)$, consider
a pair of highest weights of the form,
\[ \mu=((n-1)a, (n-2)a, \cdots, a, 0) \quad \text{and}\quad
\lambda=((n-1)b, (n-2)b, \cdots, b, 0),\]
for some positive integers $k, ~a, ~b$. Then the characters can be
expanded as Vandermonde determinants and we have, 
\[\begin{split} S_{\mu} &= \prod_{i<j}(x_i^{a+1}-x_j^{a+1})\\
\text{and} \quad S_{\lambda}& =\prod_{i<j}(x_i^{b+1}-x_j^{b+1}).
\end{split}
\]
Thus we see that $S_{\mu}$ divides $S_{\lambda}$ if $(a+1)|(b+1)$. 

It would be of interest to give necessary and sufficient criteria on the
highest weights $\mu$ and $\lambda$ to ensure that $S_{\mu}$ divides
$S_{\lambda}$.  
\end{remark}

\section{An arithmetical application} We present here an arithmetical
application to recovering l-adic representations. Corollary
\ref{adjoint} was motivated by the question of knowing the
relationship between two $l$-adic representations given that their
adjoint representations are isomorphic. On the other hand, the
application to generalised Asai representations given below, was
suggested by the work of D. Ramakrishnan. We refer to \cite{R} for
more details.  

Let $K$ be a global field and let $G_K$ denote the Galois group over
$K$ of an algebraic closure $\bar{K}$ of $K$. Let $F$ be a
non-archimedean local field of characteristic zero. Suppose 
 $$\rho_i:G_K\rightarrow GL_n(F), ~~i=1,2$$
are continuous, semisimple  representations 
of the Galois group $G_K$ into $GL_n(F)$, unramified outside a finite
set $S$ of places containing the archimedean places of $K$. Given
$\rho$, let $\chi_{\rho}$ denote the character of $\rho$. For each
finite place $v$ of $K$, we choose a place $\bar{v}$ of $\bar{K}$
dividing  $v$, and let $\sigma_{\bar{v}}\in G_K$ the corresponding
Frobenius element. If $v$ is unramified, then the value
$\chi_{\rho}(\sigma_{\bar{v}})$ depends only on $v$ and not on the
choice of $\bar{v}$, and we will denote this value by
$\chi_{\rho}(\sigma_v)$.  

Given a $l$-adic 
representation $\rho$ of $G_K$, we can construct other naturally
associated $l$-adic representations. We consider here two such
constructions: the first one, is given by the adjoint representation
\[Ad(\rho)=\rho\otimes \rho^*: G_K\to GL_{n^2}(F),\]
 where $\rho^*$ denotes the
contragredient representation of $\rho$. 

The second construction is a generalisation of Asai
representations.  Let $K/k$ be a Galois extension with Galois group
$G(K/k)$. Given $\rho$, we can associate the {\em pre-Asai representation}, 
\[As(\rho)=\otimes_{g\in G(K/k)}\rho^g,\]
where $\rho^g(\sigma)=\rho(\tilde{g}\sigma\tilde{g}^{-1}), ~\sigma \in
G_K$, and where $\tilde{g}\in G_k$ is a lift of $g\in G(K/k)$. At an
unramified place $v$ of $K$, which is split completely over a place
$u$ of  $k$, the
Asai character is given by,
\[ \chi_{As(\rho)}(\sigma_v)=\prod_{v|u}\chi_{\rho}(\sigma_v).\]
Hence we get that upto isomorphism, $As(\rho)$ does not
depend on the choice of the lifts $\tilde{g}$. If further $As(\rho)$
is irreducible, and $K/k$ ic cyclic, then $As(\rho)$ extends to a
representation of $G_k$ (called the Asai representation associated to
$\rho$ when $n=2$ and $K/k$ is quadratic).     

\begin{theorem}  Let  
 $$\rho_i:G_K\rightarrow GL_n(F), ~~i=1,2$$
be  continuous, irreducible  representations 
of the Galois group $G_K$ into $GL_n(F)$.
Let $R$ be the representation $Ad(\rho_i)$ (adjoint case)  or
$As(\rho_i)$ (Asai case) associated
to $\rho_i, ~i=1,2$. 
 
Suppose that the set of places $v$ of $K$ not in $S$, where   
\[{\rm Tr}({R\circ\rho_1}(\sigma_v))= {\rm Tr}({R\circ
\rho_2}(\sigma_v)),\]
is a set of places of positive density. Assume further that the
algebraic envelope of the image of $\rho_1$ and $\rho_2$ is connected
 and that the derived group is absolutely almost simple. 
Then the following holds:
\begin{enumerate}
\item(Adjoint case) There is a character
$\chi:G_K\to F^*$ such that
$\rho_2$ is isomorphic to either $\rho_1\otimes \chi$ or to 
$\rho_1^*\otimes \chi$.

\item(Asai case) There is a character
$\chi:G_K\to F^*$, and an element $g\in G(K/k)$ such that
$\rho_2$ is isomorphic to $\rho_1^g\otimes \chi$.
\end{enumerate}
\end{theorem}

\end{document}